\documentclass{article}[4pt]
\usepackage{amsmath, amssymb, amsthm, enumerate}
\title{A note on the Voronoi congruences and the residue of the Fermat quotient}
\author{Claire Levaillant}

\newcommand{\mpt}{\;\text{mod}\,p^3}
\newcommand{\mpd}{\;\text{mod}\,p^2}
\newcommand{\mpu}{\;\text{mod}\,p}

\newcommand{\mb}{\mathcal{B}}

\newcommand{\mpq}{\;\;\text{mod}\,p^4}

\begin{document}
\maketitle
\begin{center}Abstract\end{center} We prove a congruence on the residue of the Fermat quotient in base $a$ which arises from a generalization of the Voronoi congruences and from some other congruences on sums and weighted sums of divided Bernoulli numbers. As an application in the base $2$ case, we retrieve a congruence for the generalized harmonic number $H_{2,\frac{p-1}{2}}$, a generalization originally due to Sun of a classical congruence known since long for the harmonic number $H_{\frac{p-1}{2}}$ as a special case of the Lerch formula. We find a sharpening of the Voronoi congruences that is different from the one of Johnson and more computationally efficient. We prove an additional related congruence, which specialized to base $2$, allows to retrieve several congruences that were originally due to Lehmer. 

\section{Past interest in computing the residue of the Fermat quotient in base $a$}
Computing the residue of the Fermat quotient $q_a$ in general base $a$ has become important to mathematicians since the beginning of the twentieth century ever since young Arthur Wieferich showed in $1909$ that if the first case of Fermat's last theorem (FLT) is false for a prime $p\geq 5$ then this prime is such that $q_2=0\mpu$ \cite{WIE}. Such a prime got later called a Wieferich prime in his honor. \\

\textbf{FLT}. For $n\geq 3$, the equation $X^n+Y^n=Z^n$ has no solutions in integers X,Y,Z with $XYZ\neq 0$. \\

Since it was proven by Fermat for $n=4$ and by Euler for $n=3$, it then sufficed to prove it when $n$ is a prime $p\geq 5$. \\

\textbf{First case of FLT}. For fixed prime $p\geq 5$, there is no integer solution to $X^p+Y^p=Z^p$ with $XYZ$ prime to $p$. \\

The first Wieferich prime $1093$ was found by Meissner in $1913$ \cite{MEI} and the second $3511$ was found by Beeger in $1922$ \cite{BE}. As of today, there are no other known Wieferich primes. Wells Johson noticed in $1977$ \cite{WJ2} that the two numbers that are one less than the two known Wieferich primes have repetitions in their representations in base $2$ and other bases, such as
$$\begin{array}{ccccccccc}1092&=&10001000100&\text{in base $2$}&&3510&=&110110110110&\text{in base $2$}\\
1092&=&444&\text{in base $16$}&&3510&=&6666&\text{in base $8$}
\end{array}$$
Wieferich's result from $1909$ got extended to the other primes $q_3$ in $1910$ by Mirimanoff \cite{MIR1} (two known base $3$ Wieferich primes, namely $11$ and $1006003$ found by Kloss in $1965$), $q_5$ in $1914$ by Vandiver \cite{VA0} (six known base $5$ Wieferich primes), $q_p$ with $p\leq 31$ prime in $1917$ by Pollaczek, $q_p$ with $p\leq 89$ prime in $1988$ by Granville and Monagan \cite{GM} and finally $q_p$ with $p\leq 113$ prime in $1994$ by Suzuki \cite{SUZ}, at which point FLT was proven and interest dropped.

\section{Main theorem on the residue of the Fermat quotient in base $a$}
Let $p$ be an odd prime and let $a$ be an integer with $1\leq a\leq p-1$. The main purpose of this note is to use the results of \cite{LEV2} and \cite{LEV5} providing the respective residues of a sum of the first $(p-2)$ divided Bernoulli numbers and of a powers of $a$ weighted sum of the first $(p-2)$ divided Bernoulli numbers in order to find a nice expression for the residue of the Fermat quotient $q_a$ in base $a$. We recall these results below. In what follows, $w_p$ denotes the Wilson quotient.
\newtheorem{Result}{Result}
\begin{Result} Issued from Congruence $(9)$ in Theorem $6$ of \cite{LEV2}
$$\sum_{i=1}^{p-2}\mb_i=w_p\mpu$$
\end{Result}
\begin{Result} Issued from Proposition $2$ of \cite{LEV5}. \\Let $a$ be any integer with $1\leq a\leq p-1$. Then, we have:
$$\sum_{i=1}^{p-2}\frac{\mb_i}{a^i}=q_a+w_p\mpu$$
\end{Result}
The proof of the first fact is based on the combinatorial interpretation of the unsigned Stirling number of the first kind $\left[\begin{array}{l}p\\s\end{array}\right]$. Such a number namely counts the number of permutations of $Sym(p)$ which decompose into a product of $s$ disjoint cycles. By summing these numbers over the cycles we get the order of the symmetric group $Sym(p)$, that is $p!$. In $1900$, Glaisher provided the expansion of these numbers to the modulus $p^2$ \cite{GL}. By using Glaisher's formulas, we easily derive Result $1$.\\
The proof of the latter fact is based on the more common definition of the unsigned Stirling numbers of the first kind, expanding the falling factorial and specializing $x=-a$ in both factored and expanded form. The result then follows again from the Glaisher formulas.

The current paper arose from an attempt to find a similar congruence when the divided Bernoulli numbers are replaced with the ordinary ones.
A natural idea consisted in remembering that the residue of $B_i$ is the first potentially nonzero residue in the $p$-adic expansion of a sum of $i$-th powers of the first $(p-1)$ integers. The most general congruence for the sums of powers of integers gets provided by Zhi-Hong Sun in \cite{SU2}. In \cite{LEV4}, using Sun's method, we pushed his expansion to the modulus $p^4$:
\begin{equation}S_t=\sum_{a=1}^{p-1}a^t=pB_t+\frac{p^2}{2}tB_{t-1}+\frac{p^3}{6}t(t-1)B_{t-2}+\frac{p^4}{24}t(t-1)(t-2)B_{t-3}\;\mpq\end{equation}
This is also a good time to recall the fundamental consequence of the Von Staudt-Clausen's theorem: the denominator of a Bernoulli number $B_i$ consists of products of primes $p$ with multiplicities one, such that $p-1$ divides $i$ \cite{VS}\cite{CL}.\\\\
\textbf{In what follows, $i$ is a positive integer prime to $p$ and incongruent to $1$ modulo $(p-1)$.}\\\\ For a $p$-adic integer $x\in\mathbb{Z}_p$, we denote by $(x)_i$ its $(i+1)$-th residue in its $p$-adic Hensel expansion. See for instance \cite{GO}. We have under the given assumptions on $i$ and using the notations from before:
\begin{equation}S_i=p\,B_i\mpd\end{equation} Then,
\begin{equation}p\frac{B_i}{a^i}=\sum_{b=1}^{p-1}\left(\begin{array}{l}\frac{b}{a}\end{array}\right)^i=
\sum_{b=1}^{p-1}\left(\begin{array}{l}\frac{b}{a}\end{array}\right)_0^i+pi\sum_{b=1}^{p-1}\left(\begin{array}{l}\frac{b}{a}\end{array}\right)_1\left(\begin{array}{l}\frac{b}{a}\end{array}\right)_0^{i-1}\mpd\end{equation}
In the middle sum of $(3)$, the product $ba^{-1}$ with $a^{-1}$ the inverse of $a$ modulo $p$ must be treated modulo $p^2$, while in the first sum of the right hand side this product is now treated modulo $p$. Then, for fixed $a$, the product $ba^{-1}$ takes all the values between $1$ and $p-1$ exactly once when $b$ varies between $1$ and $p-1$. Whence, that sum is nothing else than a sum of $i$-th powers of the first $p-1$ integers which is congruent to $pB_i$ modulo $p^2$.
Moreover, the second coefficient in the $p$-adic expansion of $ba^{-1}$ can be written in terms of integer part as
\begin{equation}\left(\begin{array}{l}\frac{b}{a}\end{array}\right)_1=\left[\begin{array}{l}\frac{ba^{-1}}{p}\end{array}\right],\end{equation} where the inverse of $a$ is taken modulo $p$. Simplyfing the congruence by $p$ and dividing by $i$ leads to
\begin{equation}
\frac{\mb_i}{a^i}-\mb_i=\sum_{b=1}^{p-1}\left[\begin{array}{l}\frac{ba^{-1}}{p}\end{array}\right](ba^{-1})_0^{i-1}\mpu
\end{equation}
In the case when $i$ is even and with some slightly different assumption on $i$, namely that $i$ is incongruent to $0$ modulo $(p-1)$, we obtain the fundamental Voronoi congruence where the integer $a$ has been replaced with its inverse. This famous congruence got generalized to the modulus $p^2$ by Johnson in \cite{WJ} still for the even $i$'s and under the assumption that $i$ is incongruent to $2$ modulo $(p-1)$ and $a\neq 1$.
Formerly in \cite{VA}, Vandiver had built upon the Voronoi congruences by showing under the same assumption as Voronoi's that
\begin{equation}(1-a^i)\mb_i=\sum_{n=1}^{a}\sum_{b=1}^{\left[\begin{array}{l}\frac{np}{a}\end{array}\right]}(ab)^{i-1}\mpu\end{equation}
This is the congruence that was extensively used by many authors to compile tables for the irregular primes. K\"ummer had shown that FLT holds true when the exponent is a regular prime \cite{KUM}.\\

We now operate on Congruence $(5)$ by summing over the $i$'s for the range $2\leq i\leq p-3$. The left hand side gets processed through Results $1$ and $2$. The right hand side gets processed through the summation of a geometric series when $b$ is distinct from $a$ and directly for $b=a$. In the first case, it yields:
$$-1-(ba^{-1})_0^{-1}-(ba^{-1})_0^{-2}\mpu$$
In the second case, it yields $p-4$ that is $-4$ modulo $p$. Hence the following statement.
\newtheorem{Theorem}{Theorem}
\begin{Theorem} Let $a$ be an integer with $1\leq a\leq p-1$. The residue of the Fermat quotient in base $a$ can be computed as follows.
\begin{equation}\begin{split}
q_a=\frac{1}{2}(1-a^{-1})-\sum_{b=1}^{p-1}&\left[\begin{array}{l}\frac{b(a^{-1})_0}{p}\end{array}\right]
-a\sum_{b=1}^{p-1}\left[\begin{array}{l}\frac{b(a^{-1})_0}{p}\end{array}\right]\frac{1}{b}
\\&-a^2\sum_{b=1}^{p-1}\left[\begin{array}{l}\frac{b(a^{-1})_0}{p}\end{array}\right]\frac{1}{b^2}-\left[\begin{array}{l}\frac{a(a^{-1})_0}{p}\end{array}\right]\mpu\end{split}\end{equation}
\end{Theorem}
Our proof of the Voronoi congruence leads to a generalization modulo $p^2$ which has a different form than the one of Johnson (see Theorem $8$ of \cite{WJ}). We take the same assumptions on $i$ as before, except when $i$ is even (resp odd), we add the assumption that $i$ is incongruent to $2$ (resp $3$) modulo $(p-1)$. Our result is the following.
\begin{Theorem}
Let $i$ be a positive integer that is prime to $p$.  \\
(i) Suppose $i$ is even and $i$ is incongruent to $2$ modulo $(p-1)$. Then, we have:
\begin{equation}(a^i-1)\mb_i=\sum_{b=1}^{p-1}\Bigg[\frac{ab}{p}\Bigg](ab)_0^{i-1}+p\,\frac{i-1}{2}a^{i-2}\sum_{b=1}^{p-1}\Bigg[\frac{ab}{p}\Bigg]^2b^{i-2}\mpd\end{equation}
(ii) Suppose $i$ is odd and $i$ is incongruent to $1$ or $3$ modulo $(p-1)$. Then, we have:
\begin{equation}\begin{split}
B_{i-1}(a^i-1)=2\frac{\sum_{b=1}^{p-1}\Big[\frac{ab}{p}\Big](ab)_0^{i-1}}{p}&+(i-1)\sum_{b=1}^{p-1}\Bigg[\frac{ab}{p}\Bigg]^2(ab)_0^{i-2}\\
&+p\,\frac{(i-1)(i-2)}{3}\sum_{b=1}^{p-1}\Bigg[\frac{ab}{p}\Bigg]^3(ab)_0^{i-3}\;\mpd
\end{split}\end{equation}
\end{Theorem}
\noindent \textsc{Proof.} It is simply a matter of working modulo $p^3$ in the case when $i$ even and modulo $p^4$ in the case when $i$ is odd, after noticing that $\Bigg(\frac{b}{a}\Bigg)_2=0$. It yields $(8)$ in the case when $i$ is even and
\begin{equation}\begin{split}\frac{p}{2}B_{i-1}(a^i-1)=\sum_{b=1}^{p-1}\Bigg[\frac{ab}{p}\Bigg](ab)_0^{i-1}&
+\frac{p}{2}(i-1)\sum_{b=1}^{p-1}\Bigg[\frac{ab}{p}\Bigg]^2(ab)_0^{i-2}\\&+p^2\frac{(i-1)(i-2)}{6}
\sum_{b=1}^{p-1}\Bigg[\frac{ab}{p}\Bigg]^3(ab)_0^{i-3}\mpt\end{split}\end{equation}
in the case when $i$ is odd. The first member of the right hand side of $(10)$ is divisible by $p$ by $(5)$ since we have assumed that $i$ is incongruent to $1$ modulo $(p-1)$. We obtain $(9)$.

Instead of $(8)$, Johnson's congruence reads:
\begin{equation}
(a^i-1)\mb_i=\sum_{b=1}^{p-1}\Bigg[\frac{ab}{p}\Bigg](ab)^{i-1}-p\,\frac{i-1}{2}a^{i-2}\sum_{b=1}^{p-1}\Bigg[\frac{ab}{p}\Bigg]^2b^{i-2}\mpd
\end{equation}
Johnson's sharpening of the Voronoi congruence uses the Teichm\"uller characters while our proof is only based on congruences concerning sums of powers of the first $(p-1)$ integers. \\
Up to the sign, the second term of the right hand side of $(11)$ is identical to the one of $(8)$. By comparing both congruences, we derive an additional statement.
\newtheorem{Corollary}{Corollary}
\begin{Corollary} Let $i$ be a positive even integer that is prime to $p$ and incongruent to $2$ modulo $(p-1)$. Let $a$ be an integer with $2\leq a\leq p-1$. Then we have:
\begin{equation}
2(a^i-1)\mb_i=\sum_{b=1}^{p-1}\Big((ab)^{i-1}+(ab)_0^{i-1}\Big)\Bigg[\frac{ab}{p}\Bigg]\mpd
\end{equation}
\end{Corollary}

The next part is concerned with base $2$.

\section{Specificity of base $2$ and related developments}
A lot more studies were made in base $2$. There even exists a combinatorial interpretation for the residue of the Fermat quotient in base $2$. This residue relates to the number of permutations of the symmetric group $Sym(p-2)$ with an even number of ascents, denoted for convenience by $N_{p-2}$. It is shown in \cite{LEV2} that
\begin{equation}
q_2=(2N_{p-2})_0-1\;\mpu
\end{equation}
Thus, $p$ is a Wieferich prime if and only if the residue of twice the number of permutations of $Sym(p-2)$ with an even number of ascents is $1$.

In this part we prove two statements that are both in connection to the general case discussed in $\S\,2$. The first statement is directly linked to Theorem $1$ and concerns the residue of a sum of squared powers of the first $\frac{p-1}{2}$ reciprocals.
We provide a different proof than Zhi-Hong Sun's that the residue of the first generalized harmonic number of order $\frac{p-1}{2}$ is zero.
\newtheorem{Proposition}{Proposition}
\begin{Proposition}
$$H_{2,\frac{p-1}{2}}=\sum_{k=1}^{\frac{p-1}{2}}\frac{1}{k^2}=0\mpu$$
\end{Proposition}
It had been known since Wolstenholme \cite{WO} that when the order of the sum is rather $(p-1)$, this residue is zero, a result which got later generalized by Bayat to all the other generalized harmonic numbers of that order \cite{BA}. Sun's result listed as his Corollary $5.2$ in \cite{SU2} is much more general than the single case described above, as it deals with the other generalized harmonic numbers of order $\frac{p-1}{2}$ as well, including the odd powers and also working modulo $p^2$. Contrary to what happens with the even powers, when the powers are odd, the considered residue is not zero, starting with the harmonic number of order $\frac{p-1}{2}$ whose study goes back to the work of Eisenstein from $1850$. Eisenstein relates modulo $p$ the Fermat quotient in base $2$ with the alternating harmonic sum of order $p-1$ \cite{EI}:
\begin{equation*}
q_2=\frac{1}{2}\sum_{k=1}^{p-1}\frac{(-1)^{k-1}}{k}\;\;\mpu
\end{equation*}
From Eisenstein's formula and Wolstenholme's theorem we easily derive:
\begin{equation}H_{\frac{p-1}{2}}=\sum_{k=1}^{\frac{p-1}{2}}\frac{1}{k}=-2q_2\mpu\end{equation}
If we denote by $H^{'}_{p-1}$ the harmonic sum of order $p-1$ taken only on the odd integers, we have by Eisenstein's formula and Wolstendholme's theorem: \begin{equation*}
H^{'}_{p-1}=q_2\mpu
\end{equation*}
We note that Glaisher \cite{GL2} and Sun \cite{SU2} a century later successively extended the expansion for $H^{'}_{p-1}$ to the respective moduli $p^2$ and $p^3$. Sun's congruence reads:
\begin{equation}
H^{'}_{p-1}=q_2-\frac{1}{2}pq_2^2+\frac{1}{3}p^2q_2^3-\frac{1}{24}p^2B_{p-3}\;\mpt
\end{equation}
Congruence $(14)$ is a special case of Lerch's formula dating from $1905$ which asserts, written in the Vandiver form, that:
$$-aq_a=\sum_{u=1}^{a-1}\sum_{k=1}^{\left[\begin{array}{l}\frac{up}{a}\end{array}\right]}\frac{1}{k}\;\mpu$$
A proof of the latter formula appears for instance in \cite{VA2}.\\

The proof of Proposition $1$ goes as follows. We apply Theorem $1$ with $a=(2^{-1})_0$. We recall from \cite{LEV4} that
$$(2^{-1})_0=\frac{p+1}{2}$$
This is the first expansion of Lemma $2$ of \cite{LEV4}.
Then,
$$\Bigg[\frac{2.(2^{-1})_0}{p}\Bigg]=\Bigg[\frac{p+1}{p}\Bigg]=1$$
Moreover, if $b\leq \frac{p-1}{2}$, then $\Bigg[\frac{2b}{p}\Bigg]=0$. The other $b$'s to the exception of $(2^{-1})_0=\frac{p+1}{2}$ may be written as $p-k$ with $k$ varying from $1$ to $\frac{p-3}{2}$. For those $k$'s, we have $3\leq p-2k\leq p-2$. Then,
$$\Bigg[\frac{2(p-k)}{p}\Bigg]=1$$
It follows that
\begin{equation}
q_{(2^{-1})_0}=-\frac{1}{2}-4-\frac{p-3}{2}+2^{-1}\sum_{k=1}^{\frac{p-3}{2}}\frac{1}{k}-(2^{-1})^2\sum_{k=1}^{\frac{p-3}{2}}\frac{1}{k^2}\mpu,
\end{equation}
which can be rewritten using Congruence $(14)$ as:
\begin{equation}
q_{(2^{-1})_0}+q_2=-1-(2^{-1})^2H_{2,\frac{p-1}{2}}\mpu
\end{equation}
Further, by a classical identity on Fermat quotients, we have:
$$q_{(2^{-1})_0}+q_2=q_{2.(2^{-1})_0}=q_{p+1}=-1\mpu$$
The result of Proposition $1$ follows.
\newtheorem{Remark}{Remark}
\begin{Remark}
Doing $a=(2^{-1})_0$ in Congruence $(5)$ leads to
\begin{eqnarray}
2^i\mb_i-\mb_i&=&\sum_{b=1}^{p-1}\Bigg[\frac{2b}{p}\Bigg](2b)^{i-1}\mpu\\
&=&\sum_{k=1}^{\frac{p-1}{2}}2^{i-1}(p-k)^{i-1}\mpu\\
&=&(-1)^{i-1}\sum_{\begin{array}{l}k=2\\\text{$k$ is even}\end{array}}^{p-1}k^{i-1}\mpu\\
&=&(-1)^i\sum_{\begin{array}{l}k=1\\\text{$k$ is odd}\end{array}}^{p-2}k^{i-1}\mpu
\end{eqnarray}
In particular, we have for every positive integer $i$ even or odd, prime to $p$ and incongruent to $1$ modulo $(p-1)$:
\begin{equation}
1^{i-1}+3^{i-1}+\dots+(p-2)^{i-1}=(2^i-1)\mb_i\;\mpu
\end{equation}
\end{Remark}
\begin{Remark}
There also exists a congruence for the odd reciprocals. Indeed, by Corollary $5.2$ of \cite{SU2}, we know that ($p>5$):
\begin{equation}
\sum_{k=1}^{\frac{p-1}{2}}\frac{1}{k^i}=\begin{cases} (2^i-2)\,\mb_{p-i}\;\mpu&\text{if $i\in\lbrace 3,5,\dots, p-4\rbrace$}\\
0\qquad\qquad\;\;\;\;\mpu&\text{if $i\in\lbrace 2,4,\dots,p-5\rbrace$}\end{cases}
\end{equation}
Then, we have:
\hspace{-0.3cm}\begin{equation}
\frac{1}{1^{i-1}}+\frac{1}{3^{i-1}}+\frac{1}{5^{i-1}}+\dots+\frac{1}{(p-2)^{i-1}}=
\begin{cases}
q_2\qquad\qquad\qquad\;\;\;\mpu&\text{if $i=2$}\\
\Big(\frac{1}{2^{i-2}}-1\Big)\mb_{p+1-i}\!\,\mpu&\text{if $i=4,6,\dots,p-3$}\\
0\qquad\qquad\qquad\;\;\;\;\mpu &\text{if $i=3,5,\dots,p-4$}
\end{cases}
\end{equation}
\end{Remark}
Congruence $(20)$ has a refinement modulo $p^2$ when $i$ is even. D. Mirimanoff has shown in \cite{MIR2} that the same congruence holds modulo $p^2$ under the additional condition that $i$ is incongruent to $2$ modulo $(p-1)$. Independently, E. Lehmer showed in \cite{LEM} a congruence modulo $p^2$ for the odd $i$'s. We gather both of their results below, with a minor change of indices $k=i-1$ with respect to their respective original statements. In \cite{LEV2} we gave a common proof for both congruences. This proof is independent from the Voronoi type congruences and uses the Bernoulli polynomials.
\begin{Result} Lehmer--Mirimanoff congruences. Let $i$ be prime to $p$ and $p-1\not|i-2$. Then,
\begin{equation}
\sum_{r=1}^{\frac{p-1}{2}}r^{i-1}=\begin{cases}
(2^{-i+2}-1)B_{i-1}\frac{p}{2}\;\mpd&\text{if $i$ is odd}\;\; (I)_i\\
&\\
\big(\frac{1}{2^i}-1\big)\frac{2B_{i}}{i}\;\,\qquad\mpd&\text{if $i$ is even}\;\; (II)_i
\end{cases}
\end{equation}
\end{Result}

\begin{Remark}
The fact that Congruence $(20)$ holds modulo $p^2$ under the assumptions on $i$ that $i$ is even, both $i$ and $i-1$ are prime to $p$ and $i$ incongruent to $2$ modulo $(p-1)$
can be seen from Theorem $2$ $(i)$ applied in base $2$. Indeed, when $a=2$, Congruence $(8)$ reads:
\begin{equation}
(2^i-1)\mb_i=\sum_{k=1}^{\frac{p-1}{2}}(p-2k)^{i-1}+p\frac{i-1}{2}2^{i-2}\sum_{k=1}^{\frac{p-1}{2}}(p-k)^{i-2}\mpd
\end{equation}
Further, since $i-2$ is even, $(20)$ itself imposes that the second sum to the right hand side of $(26)$ is divisible by $p$ (applying $(20)$ with odd $i-1$ is licit since $i-1\neq 1\,\text{mod}\,(p-1)$ and $i\neq 1\mpu$ by assumption on $i$). The latter fact may also be used inside the first sum. Therefore, $(20)$ also holds modulo $p^2$ when $i$ is even, $i$ and $i-1$ are prime to $p$ and $i$ is incongruent to $2$ modulo $(p-1)$, that is Mirimanoff's congruence $(II)_i$ holds.
\end{Remark}
\noindent Lehmer's congruence $(I)_i$ can also be deduced from a weaker form of Theorem $2$ $(ii)$ in some cases as well, namely $i$ is odd, both $i$ and $i-1$ are relatively prime to $p$ and $i$ is incongruent to $1$ modulo $(p-1)$ (Congruence $(10)$ is then only needed modulo $p^2$). Further, we show that $(I)_i$ has an identical refinement modulo $p^3$ under the conditions on $i$ expressed below.
\begin{Proposition}
Let $i$ be an odd integer with $i$ incongruent to $0$ or $1$ or $2$ modulo $p$ and $i$ incongruent to $1$ or $3$ modulo $(p-1)$. \\
Then the Lehmer congruence $(I)_i$ holds modulo $p^3$.
\end{Proposition}

\noindent \textsc{Proof.} An application of Congruence $(9)$ with $a=2$ under the conditions of application on $i$, that is $i$ is odd, $i$ is prime to $p$ and incongruent to $1$ and $3$ modulo $(p-1)$, yields:
\begin{equation}
B_{i-1}(2^i-1)=2\frac{\sum_{k=1}^{\frac{p-1}{2}}(p-2k)^{i-1}}{p}-(i-1)2^{i-2}\sum_{k=1}^{\frac{p-1}{2}}k^{i-2}\mpd
\end{equation}
We omitted to write the last sum of $(9)$ as it is divisible by $p$ by the original congruence $(20)$ applied with odd $(i-2)$ whose application is licit since $i-2\neq 1\,\text{mod}\,(p-1)$ and $i\neq 2\mpu$.\\
Expanding the first sum to the right hand side of $(27)$ now yields:
\begin{equation}
\frac{\sum_{k=1}^{\frac{p-1}{2}}k^{i-1}}{p}=B_{i-1}\Big(1-\frac{1}{2^i}\Big)+\frac{3(i-1)}{4}\sum_{k=1}^{\frac{p-1}{2}}k^{i-2}\mpd
\end{equation}
By assumption, $i-1$ is even, $i\neq 1\mpu$ and $i-1\neq 2\,\text{mod}\,(p-1)$. Then, $(II)_{i-1}$ applies. It yields:
\begin{equation}
\sum_{k=1}^{\frac{p-1}{2}}k^{i-2}=\Big(\frac{1}{2^{i-1}}-1\Big)2\mb_{i-1}\mpd
\end{equation}
The result then follows from gathering Congruences $(28)$ and $(29)$. \hfill $\square$

\begin{Remark}
The method of \cite{LEV2} allows to generalize such congruences modulo the other $p^r$ with $r\geq 4$ as well.
\end{Remark}
We now apply again Theorem $2$ $(ii)$ under its conditions of application on $i$ and further impose that both $i-1$ and $i-2$ are relatively prime to $p$. It comes:
\begin{equation}
pB_{i-1}2^i=pB_{i-1}+2\sum_{k=1}^{\frac{p-1}{2}}(p-2k)^{i-1}-p(i-1)2^{i-2}\sum_{k=1}^{\frac{p-1}{2}}k^{i-2}\mpt
\end{equation}

\noindent Like before, the last sum of $(9)$ is divisible by $p$ and thus vanishes from the congruence. This time, contrary to before, we do not expand the first sum. But we treat the second sum just like before, using Congruence $(29)$. After regrouping the different terms we obtain:
\begin{equation}
pB_{i-1}=\frac{1}{2^{i-2}}\sum_{k=1}^{\frac{p-1}{2}}(p-2k)^{i-1}\mpt
\end{equation}
This congruence was originally proven by Emma Lehmer in \cite{LEM}. Her assumptions are weaker than ours: she only assumes that $i$ is odd and incongruent to $3$ modulo $(p-1)$.\\

\textsc{Email:} \textit{clairelevaillant@yahoo.fr}

\end{document}